\documentclass[11pt]{amsart}
\usepackage{amsmath,amssymb,color,cite,placeins,multirow,mathtools,colonequals,verbatim}
\usepackage{etoolbox}
\apptocmd{\sloppy}{\hbadness 10000\relax}{}{}

\makeatletter
\@namedef{subjclassname@2020}{\textup{2020} Mathematics Subject Classification}
\makeatother

\numberwithin{equation}{section}

\newtheorem{thm}[equation]{Theorem}
\newtheorem{prop}[equation]{Proposition}
\newtheorem{lemma}[equation]{Lemma}
\newtheorem{cor}[equation]{Corollary}

\theoremstyle{definition}
\newtheorem{rmk}[equation]{Remark}

\newcommand{\F}{\mathbb{F}}

\newcommand{\Z}{\mathbb{Z}}

\DeclareMathOperator{\charp}{char}

\DeclareMathOperator{\Tr}{Tr}

\newcommand{\mybar}[1]{#1\llap{$\overline{\phantom{\rm#1}}$}}
\newcommand{\abs}[1]{\lvert #1 \rvert}

\usepackage[colorlinks,pagebackref,pdftex, bookmarks=false]{hyperref}

\begin{document}

\author{Zhiguo Ding}
\address{
  Department of Mathematics,
  University of Michigan,
  530 Church Street,
  Ann Arbor, MI 48109-1043 USA
}
\email{dingz@umich.edu}

\author{Michael E. Zieve}
\address{
  Department of Mathematics,
  University of Michigan,
  530 Church Street,
  Ann Arbor, MI 48109-1043 USA
}
\email{zieve@umich.edu}

\title[Sums of powers of consecutive squares]
{On sums of powers of consecutive squares over finite fields, and sums of distinct values of polynomials}

%\date{\today}

\begin{abstract}
For each odd prime power $q$, and each integer $k$, we determine the sum of the $k$-th powers of all elements $x\in\F_q$ for which both $x$ and $x+1$ are squares in $\F_q^*$.  We also solve the analogous problem when one or both of $x$ and $x+1$ is a nonsquare.  We use these results to determine the sum of the elements of the image set $f(\F_q)$ for each $f(X)\in\F_q[X]$ of the form $X^4+aX^2+b$, which resolves two conjectures by Finch-Smith, Harrington, and Wong.
\end{abstract}


\maketitle

%#######################################################################
%#######################################################################

\section{Introduction}

Let $p$ be an odd prime.  The problem of counting or estimating the number of tuples of $m$ consecutive squares in $(\Z/p\Z)^*$ has a rich history, with contributions from Aladov \cite{A}, Jacobsthal \cite{J}, Davenport \cite{D1,D2,D3}, Weil \cite{W}, and others.  It is natural to seek further information about $m$-tuples of consecutive squares, beyond merely their number.
In this paper we determine the moments of the set of initial elements of such tuples when $m=2$.  That is, we determine the sum of the $k$-th powers of all elements $x\in\Z/p\Z$ for which both $x$ and $x+1$ are nonzero squares, for any prescribed integer $k$.
In fact, we solve the analogous problem where $\Z/p\Z$ is replaced by an arbitrary finite field $\F_q$ of odd order, in which case we write $\chi(x)$ for the quadratic character on $\F_q$, which is defined by $\chi(x):=x^{(q-1)/2}$ for any $x\in\F_q$.  Thus $\chi(x)=1$ if $x$ is a nonzero square in $\F_q$, $\chi(x)=-1$ if $x$ is a nonsquare in $\F_q$, and $\chi(0)=0$.

\begin{thm}\label{main}
Let $q$ be an odd prime power and pick $a,b\in\F_q^*$.  Write
\[
D\colonequals \{x\in\F_q\colon\chi(x)=\chi(a) \text{ and }\chi(x+1)=\chi(b)\}.
\]
Let $k$ be an integer not divisible by $(q-1)/2$, and let $\ell$ and $\varepsilon$ be the unique integers such that $k\equiv \varepsilon \ell \pmod{q-1}$,\, $\varepsilon\in\{1,-1\}$, and $0<\ell<(q-1)/2$.  Then
\[
\sum_{x\in D} x^k = \frac{1+\chi(-a)+ 2^{-\ell}\cdot\chi(c)\cdot\sum_{i=0}^{\lfloor \ell/2\rfloor} 4^{-i}\binom{\ell}{2i}\binom{2i}{i}}{4\cdot (-1)^{\ell+1}},
\]
where $c:=ab$ if $\varepsilon=1$ and $c:=b$ if $\varepsilon=-1$.
In particular, if $q>3$ then
\[
\sum_{x\in D}x = \frac{1+\chi(-a)}4+\frac{\chi(ab)}8,
\]
and if $q>5$ then
\[
\sum_{x\in D}x^2 = \frac{-1-\chi(-a)}4-3\cdot \frac{\chi(ab)}{32}.
\]
\end{thm}

\begin{rmk}
The analogous result when $(q-1)/2$ divides $k$ is easy.  For, if $k=m(q-1)/2$ with $m\in\Z$ then $\sum_{x\in D}x^k = \chi(a)^m\cdot\abs{D}$, and it is well-known that $4\abs{D}=q-2-\chi(-a)-\chi(b)-\chi(ab)$ (cf.\ Remark~\ref{consec}).
\end{rmk}

\begin{comment}

Magma verification for small q and for |ell| < (q-1)/2:

for q in [3..1001 by 2] do if IsPrimePower(q) then K:=GF(q); S:={i^2:i in K}; N:={i:i in K|i notin S}; S diff:={K!0}; T:={S,N}; e:=(q-1) div 2;
for k in [1..e-1] do
s:=1/(K!2)^k*&+[Binomial(k,2*i)*Binomial(2*i,i)/(K!4)^i:i in [0..Floor(k/2)]];
for eps in {1,-1} do ell:=k*eps;
for A,B in T do a:=Random(A); b:=Random(B);
D:={x:x in A|x+1 in B};
if #D eq 0 then Sf:=0; else
Sf:=&+[x^ell:x in D]; end if;
if eps eq 1 then c:=a*b; else c:=b; end if;
if Sf ne (1+(-a)^e+c^e*s)/(K!4*(-1)^(k+1)) then q,a,b,"wrong"; break q; end if; 
if ell eq 1 and Sf ne (1+(-a)^((q-1) div 2))/4+(a*b)^((q-1) div 2)/8 then q,a,b,"wrong1"; break q; end if;
if ell eq 2 and Sf ne (-1-(-a)^((q-1) div 2))/4-3*(a*b)^((q-1) div 2)/32 then q,a,b,"wrong2"; break q; end if;
end for; end for; end for; "tested",q; end if; end for;

\end{comment}

More generally, it seems natural to study the sum of the values of some simple function (such as the identity function, or powers of the identity) on some ``nice'' subset of $\F_q$ (such as the set of elements $x$ in $\F_q$ for which $\chi(x)=\chi(x+1)=1$).
The case that the function is the constant function $1$ amounts to determining the size of the set, so the classical question of determining the number of consecutive squares mod $p$, or determining the sizes of other natural subsets of $\F_p$, may be viewed as special cases of this general class of questions.

Another instance of this general class of questions dates back to work of Stern \cite{Stern}, who determined the sum of the ``triangular'' elements of $\F_p$ for any odd prime $p$, namely the elements of $\F_p$ which can be written as $(x^2+x)/2$ for some $x\in\F_p$.  
Shorter proofs of Stern's results were given in modern language by Stetson \cite{Stetson} and Gross et al.~\cite{GHM}, with the latter authors also solving the same question with $(x^2+x)/2$ replaced by any quadratic polynomial $ax^2+bx+c$.  
More generally, for any $f(X)\in\F_q[X]$, let $f(\F_q)$ be the value set $\{f(x):x\in\F_q\}$, and let $S(f)$ be the sum of the elements of $f(\F_q)$.  When $q$ is prime and $\deg(f)=3$, the quantity $S(f)$ was determined by Finch-Smith et al.~\cite{FHW}, who also proposed two conjectures about the degree-$4$ case.  We will use Theorem~\ref{main} to prove generalizations of these two conjectures.

\begin{thm}\label{conj3}
For any odd prime power $q$ with $q>5$, and any $a,b\in\F_q$ with $a\ne 0$, the quantity $S(X^4+aX^2+b)$ equals 
\[
-\Bigl( \frac{4+\chi(-1)+4\chi(-2a)}{64} \Bigr) \cdot a^2 + \Bigl( \frac{4+\chi(-1)-2\chi(-a)+2\chi(-2a)}{8} \Bigr) \cdot b.
\]
%
\begin{comment}
or equivalently, depending on one's taste,
\[
\frac{-a^2+8b}{16} + \chi(-1)\cdot\frac{-a^2+8b}{64} + \chi(-a)\cdot\frac{-b}{4} + \chi(-2a)\cdot\frac{-a^2+4b}{16}
\]
\end{comment}
%
\begin{comment}
% version which doesn't use chi:
\begin{align*}
&\frac{-9c^2+40e}{64}&&\text{if $q\equiv 1\pmod 8$ and $c$ is a square in\/ $\F_q^*$;} \\
&\frac{-c^2+40e}{64}&&\text{if $q\equiv 1\pmod 8$ and $c$ is a nonsquare in\/ $\F_q^*$;} \\
&\frac{-7c^2+56e}{64}&&\text{if $q\equiv 3\pmod 8$ and $c$ is a square in\/ $\F_q^*$;} \\
&\frac{c^2-8e}{64}&&\text{if $q\equiv 3\pmod 8$ and $c$ is a nonsquare in\/ $\F_q^*$;} \\
&\frac{-c^2+8e}{64}&&\text{if $q\equiv 5\pmod 8$ and $c$ is a square in\/ $\F_q^*$;} \\
&\frac{-9c^2+72e}{64}&&\text{if $q\equiv 5\pmod 8$ and $c$ is a nonsquare in\/ $\F_q^*$;} \\
&\frac{c^2+24e}{64}&&\text{if $q\equiv 7\pmod 8$ and $c$ is a square in\/ $\F_q^*$;} \\
&\frac{-7c^2+24e}{64}&&\text{if $q\equiv 7\pmod 8$ and $c$ is a nonsquare in\/ $\F_q^*$.}
\end{align*}
\end{comment}
\end{thm}

%
\begin{comment}

Check that the formula involving chi agrees with the eight cases of Conjecture 3 of FHW:

qs:=[9,3,5,7]; ct:=1;
for q in qs do K:=GF(q); s:=(q-1) div 2; P<c,e>:=PolynomialRing(K,2);
conj:=[(-9*c^2+40*e)/64, (-c^2+40*e)/64, (-7*c^2+56*e)/64, (c^2-8*e)/64, (-c^2+8*e)/64, (-9*c^2+72*e)/64, (c^2+24*e)/64, (-7*c^2+24*e)/64];
for c0 in [K!1,PrimitiveElement(K)] do
-(4+(-1)^s+4*(-2*c0)^s)/64*c^2 + (4+(-1)^s-2*(-c0)^s+2*(-2*c0)^s)/8*e eq conj[ct]; ct+:=1; end for; end for;

\end{comment}
%

Theorem~\ref{conj3} is a generalized and corrected version of Conjecture~3 of \cite{FHW}.
%
% Justify:
% Corrects: Conj.3 of FHW false when p=3 or p=5
% Generalizes: FHW require our prime power q to be prime.
%
For completeness, we note that $S(X^4+b)$ is determined in Proposition~\ref{power}, and also that if $q$ is even then $S(X^4+aX^2+b)=S(X^2+aX+b)$ is determined in Proposition~\ref{3}.

As a consequence of Theorem~\ref{conj3}, we determine the collection of possibilities for $S(X^4+aX^2)$ when $a$ varies over $\F_q$.

\begin{cor}\label{conj2}
For any odd prime power $q$ with $q>5$, the set
\[
T \colonequals \{ S(X^4+aX^2) : a \in \F_q\}
\]
equals
\begin{enumerate}
\item all of\/ $\F_q$ if $q\equiv 3 \pmod 4$ and $q\equiv 3, 5 \text{ or } 6 \pmod 7$;
\item the set of squares in\/ $\F_q$ if either $q\equiv 1\pmod{12}$ or both $q\equiv 3\pmod 4$ and $q\equiv 0, 1, 2 \text{ or } 4 \pmod 7$;
\item the set of fourth powers in\/ $\F_q$ if $q\equiv 5 \text{ or } 21 \pmod {24}$;
\item the set of squares in\/ $\F_q$ which are not fourth powers in $\F_q^*$ if $q\equiv 9 \text{ or } 17 \pmod {24}$.
\end{enumerate}
Moreover, if $q\equiv 3 \pmod 4$ then $S(X^4+8X^2) =1$.
\end{cor}

Corollary~\ref{conj2} is a generalized, corrected, and sharpened version of Conjecture~2 of \cite{FHW}.
%
% Justification of the above sentence:
%
% Corrects: e.g. p=17 is a counterexample to Conj.2 of FHW, since S(f_0)=0 is in R(g_4) and hence is not in R(g_2)\R(g_4).  Also S(f_8)=0\ne 1 when p=3.
% Sharpens: we replace the conditions about whether -1 is in S with explicit congruences on q modulo 3 or 7.
% Generalizes: they required our prime power q to be a prime.

For completeness, we provide a relatively simple result determining $S(f)$ when $\deg(f)\le 3$.

\begin{prop}\label{3}
For any prime power $q$ with $q>4$, and any $a,b,c,d\in\F_q$ with $a\ne 0$, we have
\begin{enumerate}
\item $S(aX+b)=0$;
\item $S(aX^2+bX+c)$ equals
\begin{itemize}
\item $(4ac-b^2)/(8a)$ if $q$ is odd;
\item $0$ if $q$ is even;
\end{itemize}
\item $S(aX^3+bX^2+cX+d)$ equals
\begin{itemize}
\item $\lfloor (2q+1)/3\rfloor\cdot (27a^2d-9abc+2b^3)/(27a^2)$ if $b^2\ne 3ac$ and $3\nmid q$;
\item $0$ if $b^2=3ac$ and $q\equiv 2\pmod 3$;
\item $2(27a^2d-b^3)/(81a^2)$ if $b^2=3ac$ and $q\equiv 1\pmod 3$;
\item $0$ if $27\mid q$;
\item $-b^3/a^2$ if $q=9$.
\end{itemize}
\end{enumerate}
%
\begin{comment}
% old version, including the values q<=4:
\begin{enumerate}
\item $S(aX+b)$ equals $0$ if $q>2$, and equals $1$ if $q=2$;
\item $S(aX^2+bX+c)$ equals
\begin{enumerate}
\item $(4ac-b^2)/(8a)$ if $q$ is odd and $q>3$;
\item $ab^2+a-c$ if $q=3$;
\item $0$ if $b\ne 0$ and $q$ is even and $q>4$;
\item $a^2b^2$ if $b\ne 0$ and $q=4$;
\item $c$ if $b\ne 0$ and $q=2$;
\item $0$ if $b=0$ and $q$ is even and $q>2$;
\item $1$ if $b=0$ and $q=2$;
\end{enumerate}
\item $S(aX^3+bX^2+cX+d)$ equals
\begin{enumerate}
\item $\lfloor (2q+1)/3\rfloor\cdot (27a^2d-9abc+2b^3)/(27a^2)$ if $b^2\ne 3ac$ and $3\nmid q$ and $q>4$;
\item $a^2bc+a+d$ if $b^2\ne ac$ and $q=4$;
\item $d$ if $b\ne c$ and $q=2$;
\item $0$ if $b^2=3ac$ and $q\equiv 2\pmod 3$ and $q>2$;
\item $2(27a^2d-b^3)/(81a^2)$ if $b^2=3ac$ and $q\equiv 1\pmod 3$ and $q>4$;
\item $a$ if $b^2=ac$ and $q\mid 4$;
\item $0$ if $27\mid q$;
\item $-b^3/a^2$ if $q=9$;
\item $-abc-b+bc^2-d$ if $b\ne 0$ and $q=3$;
\item $0$ if $b=0$ and $\chi(-c/a)\ne 1$ and $q=3$;
\item $d$ if $b=0$ and $\chi(-c/a)=1$ and $q=3$.
\end{enumerate}
\end{enumerate}
\end{comment}
\end{prop}

Proposition~\ref{3} subsumes and generalizes \cite[Thm.~2]{FHW}, \cite[Thms.~1.2 and ~3.2]{GHM}, \cite[Satz in \S 4]{Stern}, and \cite[Thm.~I]{Stetson}.  It turns out that Proposition~\ref{3} follows easily from known results when $\F_q$ has characteristic at least $5$, but further arguments are required in characteristics $2$ and $3$.

We conclude this introduction by listing some classes of polynomials $f(X)\in\F_q[X]$ of arbitrary degree for which it is easy to compute $S(f)$.

\begin{prop}\label{power}
For any prime power $q$, let $f(X)=aX^n+b$ for some positive integer $n$ and some $a,b\in\F_q$ with $a\ne 0$.  Then 
\begin{enumerate}
\item $S(f) = a+2b$ if $(q-1)\mid n$;
\item $S(f) = b \cdot \Bigl(1-\frac{1}{\gcd(n,q-1)}\Bigr)$ if $(q-1)\nmid n$.
\end{enumerate}
\end{prop}

\begin{prop}\label{odd}
For any prime power $q$, let $f(X)=X^r B(X^s)$ for some $B(X)\in\F_q[X]$, some positive divisor $s$ of $q-1$ and some positive integer $r$ with $s\nmid r$.  Then $S(f)=0$.
\end{prop}

\begin{cor}\label{oddplus}
For any $f(X)$ as in Proposition~\ref{odd}, and any $a\in\F_q$, if $g(X):=f(X)+a$ then $S(g)=a\cdot\abs{f(\F_q)}$.
\end{cor}

The expression in Corollary~\ref{oddplus} becomes explicit for any class of polynomials $f(X)$ in Proposition~\ref{odd} for which one knows an explicit formula for $\abs{f(\F_q)}$.

\begin{rmk}
Proposition~\ref{power} generalizes and corrects \cite[Thm.~3]{FHW}.  Proposition~\ref{odd} improves \cite[Thm.~5]{FHW} in multiple ways.
\end{rmk}

%
% Justify: Thm.3 of FHW is false when a=0, b=1, and r=p-1, where f(F_p)={1} so that S(f)=1, but Thm.3 of FHW says S(f)=a+2b=2\ne 1.
% Thm.5 of FHW requires q to be prime and gcd(d,r)=1, so that Prop.~\ref{odd} generalizes it even if we restrict to prime q.

This paper is organized as follows.  In the next section we prove Theorem~\ref{main}.
Then in Section~\ref{sec3} we prove Theorem~\ref{conj3} and Corollary~\ref{conj2}, and in Section~\ref{sec4} we prove Propositions~\ref{3}--\ref{odd} and Corollary~\ref{oddplus}.

%######################################################
%######################################################

\section{Proof of Theorem~\ref{main}}

\begin{proof}[Proof of Theorem~\ref{main}]
Since the last sentence in Theorem~\ref{main} consists of the instances $k=1$ and $k=2$ of the next-to-last sentence, it suffices to prove the next-to-last sentence.  Moreover, it suffices to prove the result when $k=\ell$, since $x^k=x^{\varepsilon\ell}$ for $x\in D$, and also
$\sum_{x\in D}x^{-\ell}=\sum_{x\in D}(x^{-1})^{\ell}$ and
\[
\{x^{-1}:x \in D\}=\{x\in\F_q:\chi(x)=\chi(a) \text{ and }\chi(x+1)=\chi(ab)\}.
\]
Henceforth we assume that $k=\ell$, so that $0<k<(q-1)/2$.

We now show that it suffices to treat the case $b=a$.  Pick $\zeta\in\F_q^*$ of order $(q-1)/2$.  Then $\chi(\zeta)=1$, and since $0<k<(q-1)/2$ we have $\zeta^k\ne 1$.  Since $Z:=\{x\in\F_q:\chi(x)=\chi(a)\}$ is preserved by multiplying by $\zeta$, the sum $S:=\sum_{x\in Z} x^k$ satisfies
\[
S = \sum_{x\in Z} (\zeta x)^k = \zeta^k S,
\]
so that $S=0$.  Note that $Z$ is the union of the pairwise disjoint sets $D$, $E$, and $F$, where
\[
E:=\{x\in\F_q:\chi(x)=\chi(a) \text{ and } \chi(x+1)=-\chi(b)\}
\]
and
\[
F:=\{x\in\F_q:\chi(x)=\chi(a) \text{ and }\chi(x+1)=0\}.
\]
The equality $S=0$ says that
\[
\sum_{x\in E} x^k = - \sum_{x\in D} x^k - \sum_{x\in F} x^k.
\]
Note that $F=\{-1\}$ if $\chi(-1)=\chi(a)$, and $F$ is empty if $\chi(-1)=-\chi(a)$.  Thus in any case $\sum_{x\in F} x^k = (-1)^k \bigl(1+\chi(-a)\bigr)/2$, so that
\[
\sum_{x\in E} x^k = - \sum_{x\in D} x^k - (-1)^k\frac{1+\chi(-a)}2.
\]
In light of this, Theorem~\ref{main} is true for all pairs $(a,b)$ if it is true for the special pairs $(a,a)$.

Henceforth suppose $b=a$, and write
\[
C\colonequals \{ (u,v) : u,v \in \F_q^*,\, u^2-v^2=a \}.
\]
Let $\Lambda$ be the set of elements $t\in\F_q^*$ for which $t^4=a^2$, or equivalently $t^2\in\{a,-a\}$.  For $t \in \F_q^*\setminus \Lambda$, define
\[
\theta(t) \colonequals \Bigl( \frac{t^2+a}{2t}, \frac{t^2-a}{2t} \Bigr).
\]

We now show that $\theta$ induces a bijection $\F_q^*\setminus\Lambda\to C$.  Plainly $\theta(t)\in C$, since
\[
\frac{(t^2+a)^2}{4t^2} - \frac{(t^2-a)^2}{4t^2} = a.
\]
Next, the sum of the two components of $\theta(t)$ is $t$, so that $\theta$ is injective on $\F_q^*\setminus \Lambda$. Finally, we show that $\theta$ induces a surjection from $\F_q^* \setminus \Lambda$ onto $C$: for any $(u,v) \in C$, pick $t \in \mybar\F_q$ with $t^2 - a = 2vt$.  Plainly $t\ne 0$, so that $v=(t^2-a)/(2t)$.  
Then  $u^2 = v^2 + a = ( (t^2+a)/(2t) )^2$, so that $u = \pm (t^2+a)/(2t)$.  Note that $\hat{t}\colonequals -a/t$ satisfies $(\hat{t}^2-a)/(2\hat{t}) = v$ and $(\hat{t}^2+a)/(2\hat{t}) = -(t^2+a)/(2t)$, so upon replacing $t$ by $\hat{t}$ if necessary we may assume that $v=(t^2-a)/(2t)$ and $u=(t^2+a)/(2t)$, whence $t^2\ne \pm a$ and $u+v=t$, so that $t\in\F_q^*\setminus \Lambda$ and $(u,v)=\theta(t)$.  Thus $\theta$ defines a bijection $\F_q^* \setminus \Lambda \to C$.

Next we show that $\phi\colon (u,v)\mapsto v^2/a$ defines a surjective function $C\to D$ in which each element of $D$ has exactly four preimages.  If $(u,v)\in C$ then $x\colonequals v^2/a$ satisfies $x+1=u^2/a$ and $\chi(x)=\chi(a)=\chi(x+1)$, so that $x\in D$.  Conversely, if $x\in D$ then $\phi^{-1}(x)$ consists of the pairs $(u,v)$ of elements of $\F_q^*$ such that $u^2-v^2=a$ and $v^2/a=x$.  
Since $\chi(ax)=1$, the equality $v^2/a=x$ holds for exactly two elements $v \in \F_q^*$, and for each such $v$ we have $v^2+a=ax+a=a(x+1)$ so that $\chi(v^2+a)=\chi(a(x+1))=1$, so the equality $u^2=v^2+a$ holds for exactly two elements $u \in \F_q^*$.  Thus $\abs{\phi^{-1}(x)}=4$, so that $\phi$ defines a surjective $4$-to-$1$ function $C\to D$.

\begin{comment}
This proves the claim, which implies that
\[
\abs{D} = \frac{\abs{C}}4 = \frac{q-1 - \abs{\Lambda}}4.
\]
Note that $\abs{\Lambda}=2$ if $q\equiv 3 \pmod 4$, since then exactly one of $a$ and $-a$ is a square;
and if $q\equiv 1\pmod 4$ then $\abs{\Lambda}$ equals $4$ if $a$ is a square and $0$ if $a$ is a nonsquare.
\end{comment}

Now we compute
\begin{equation}\label{fim1}
4 \sum_{x\in D} x^k = \sum_{t \in \F_q^* \setminus \Lambda} \phi(\theta(t))^k = \sum_{t\in\F_q^*\setminus \Lambda} \frac{(t^2/a-2+a/t^2)^k}{4^k}.
\end{equation}
For any integer $i$, the value of $\sum_{t\in\F_q^*} t^i$ is $-1$ if $(q-1)\mid i$ and $0$ otherwise.  Since $0<2k<q-1$, it follows that $\sum_{t\in\F_q^*} (t^2/a-2+a/t^2)^k$ equals the negative of the constant term of the Laurent polynomial
\begin{align*}
\Bigl(\frac{X^2}a-2+\frac{a}{X^2}\Bigr)^k &= \sum_{j=0}^k \binom{k}{j} (-2)^{k-j} \Bigl(\frac{X^2}a+\frac{a}{X^2}\Bigr)^j \\
&= \sum_{j=0}^k \binom{k}{j} (-2)^{k-j} \sum_{i=0}^j \binom{j}{i} a^{j-2i} X^{4i-2j} \\
&= \sum_{i=0}^k \sum_{j=i}^k \binom{k}{j} \binom{j}{i} (-2)^{k-j} a^{j-2i} X^{4i-2j}.
\end{align*}
This constant term is
\[
\sum_{i=0}^{\lfloor k/2\rfloor} \binom{k}{2i} \binom{2i}{i} (-2)^{k-2i},
\]
so that
\[
4^{k+1}\sum_{x\in D} x^k = - \sum_{i=0}^{\lfloor k/2\rfloor} \binom{k}{2i} \binom{2i}{i} (-2)^{k-2i} -\sum_{t\in\Lambda} \Bigl(\frac{t^2}a-2+\frac{a}{t^2}\Bigr)^k.
\]
Each $t\in\Lambda$ satisfies $t^2=\pm a$.  If $t^2=a$ then $t^2/a-2+a/t^2=0$, and if $t^2=-a$ then $t^2/a-2+a/t^2=-4$.  Since in addition there are $1+\chi(-a)$ elements $t\in\Lambda$ for which $t^2=-a$, we conclude that
\[
\sum_{t\in\Lambda} \Bigl(\frac{t^2}a-2+\frac{a}{t^2}\Bigr)^k = \bigl(1+\chi(-a)\bigr) (-4)^k.
\]
Thus
\[
4^{k+1}\sum_{x\in D} x^k = - \sum_{i=0}^{\lfloor k/2\rfloor} \binom{k}{2i} \binom{2i}{i} (-2)^{k-2i} - \bigl(1+\chi(-a)\bigr) (-4)^k,
\]
or equivalently
\[
\sum_{x\in D} x^k = \frac{1+\chi(-a)  + 2^{-k}\cdot \sum_{i=0}^{\lfloor k/2\rfloor} 4^{-i} \binom{k}{2i} \binom{2i}{i}}{4\cdot (-1)^{k+1}}.
\]
Since $b=a$ we have $\chi(ab)=1$, so this matches the desired equality, which concludes the proof.
\end{proof}

\begin{rmk}\label{consec}
As a by-product, the above proof yields the classical formula for $\abs{D}$ in case $b=a$: for, $\phi\circ\theta$ is a surjective $4$-to-$1$ function $\F_q^*\setminus\Lambda\to D$, so since $\abs{\Lambda}=2+\chi(a)+\chi(-a)$ we obtain $4\abs{D}=q-3-\chi(a)-\chi(-a)$.
More generally, this is also the formula for $\abs{D}$ for any $a,b\in\F_q^*$ with $\chi(b)=\chi(a)$.  We now deduce from this the formula for $\abs{D}$ when $\chi(b)=-\chi(a)$. Let $Z$ be the set of $x\in\F_q$ with $\chi(x)=\chi(a)$.
Since the number of $x\in Z$ for which $\chi(x+1)\ne 0$ is $(q-1)/2-(1+\chi(-a))/2$, then the number of $x\in Z$ for which $\chi(x+1)=-\chi(a)$ is $(q-1+\chi(a)-\chi(-a))/4$.  Thus for any $a,b\in\F_q^*$ we have $4\abs{D}=q-2-\chi(-a)-\chi(b)-\chi(ab)$.
\end{rmk}

%
\begin{comment}

Verify this via Magma for small q.

for q in [3..10^5 by 2] do if IsPrimePower(q) then K:=GF(q); S:={i^2:i in K}; N:={i:i in K|i notin S}; S diff:={K!0}; for A,B in {S,N} do a:=Random(A); b:=Random(B); D:={i:i in A|i+1 in B};
if a*b in S then chiab:=1; else chiab:=-1; end if;
if b in S then chib:=1; else chib:=-1; end if;
if -a in S then chima:=1; else chima:=-1; end if;
if 4*#D ne q-2-chiab-chib-chima then q,a,b,"WRONG"; break q; end if; end for; "tested",q; end if; end for;

\end{comment}
%

%######################################################
%######################################################

\section{Proofs of Theorem~\ref{conj3} and Corollary~\ref{conj2}}\label{sec3}

In this section we prove Theorem~\ref{conj3} and Corollary~\ref{conj2}.  We begin with the following generalization of \cite[Thm.~2.1]{Sun}.

\begin{lemma}\label{4val}
For any odd prime power $q$, pick $a\in\F_q^*$ and write $f(X):=X^4+aX^2$.  Then $N:=\abs{f(\F_q)}$ satisfies
\[
N = \frac{3q+4+\chi(-1)-2\chi(-a)+2\chi(-2a)}8 = \Bigl\lfloor\frac{3q+7-2\chi(-a)}8\Bigr\rfloor.
\]
\end{lemma}
 
\begin{proof}
Let $Z$ be the set of squares in $\F_q^*$.  Then $f(\F_q)$ consists of $f(0)=0$ and $W:=\{x^2+ax:x\in Z\}$.  Plainly $0\in W$ if and only if $-a\in Z$, so that $N=\abs{W}+(1-\chi(-a))/2$.  For $x,y\in Z$ we have $x^2+ax=y^2+ay$ if and only if $y\in\{x,-a-x\}$.  
Note that $x=-a-x$ just when $x=-a/2$, and that $-a/2\ne -a$.  Writing $M$ for the number of $x\in Z$ for which $y:=-a-x$ is also in $Z$, it follows that $\abs{W}=\abs{Z}-M/2+(1+\chi(-a/2))/4$.  
Next, writing $x:=a\tilde{x}$, we see that $M$ is the number of $\tilde{x}\in\F_q^*$ for which $\chi(\tilde{x})=\chi(a)=\chi(-1-\tilde{x})$, or equivalently $\chi(\tilde{x})=\chi(a)$ and $\chi(\tilde{x}+1)=\chi(-a)$.  Thus Remark~\ref{consec} determines $M$, which in turn determines $\abs{W}$ and $N$.
\end{proof}

%
\begin{comment}

Magma verification for small q:

for q in [3..1001 by 2] do if IsPrimePower(q) then K:=GF(q); for c in K do if c ne 0 then N:=#{x^4+c*x^2:x in K};
if IsSquare(-K!1) then s:=1; else s:=-1; end if;
if IsSquare(-c) then t:=1; else t:=-1; end if;
if IsSquare(-2*c) then u:=1; else u:=-1; end if;
if #{N,(3*q+4+s-2*t+2*u)/8,Floor((3*q+7-2*t)/8)} ne 1 then q,c,"WRONG"; break q; end if; end if; end for; "tested",q; end if; end for;

\end{comment}
%

\begin{proof}[Proof of Theorem~\ref{conj3}]
For $N$ as in Lemma~\ref{4val}, we have $S(X^4 + aX^2 + b) = S(X^4 + aX^2) + bN$.  Since Lemma~\ref{4val} determines $N$, it remains to compute $S(X^4+aX^2)$.  Note that
\[
\{x^4+ax^2 : x \in \F_q\} = \{0\} \cup \{x^2+ax : x \in \F_q, \chi(x)=1\}.
\]
Since $q>5$, we have $\sum_{\chi(x)=1} (x^2+ax) = 0$.  Since $x^2+ax=y^2+ay$ if and only if $y\in\{x,-a-x\}$, this says that
\[
S(X^4+aX^2) = -\frac12\cdot \!\!\!\!\sum_{\substack{\chi(x)=1 \\ \chi(-x-a)=1}} \!\!\!(x^2+ax) \,\,+\frac12\cdot \!\!\!\sum_{\substack{\chi(x)=1 \\ x=-x-a}} \!(x^2+ax).
\]
Writing $x=a\tilde{x}$, it follows that
\[
S(X^4+aX^2) = -\frac{a^2}2\cdot \!\!\!\!\!\!\sum_{\substack{\chi(\tilde{x})=\chi(a) \\ \chi(\tilde{x}+1)=\chi(-a)}} \!\!\!\!\!\! (\tilde{x}^2+\tilde{x}) + \frac{a^2}2\cdot \!\!\!\sum_{\substack{\chi(\tilde{x})=\chi(a) \\ \tilde{x}=-\tilde{x}-1}} \!\!\! (\tilde{x}^2+\tilde{x}).
\]
The first sum is determined in Theorem~\ref{main}, and the second sum equals $-a^2/8$ if $\chi(-2a)=1$, and $0$ otherwise, so it always equals $-a^2 (1+\chi(-2a)) / 16$.  Theorem~\ref{conj3} follows.
\end{proof}

%
\begin{comment}

If $q=3$ then the image of $x^4+cx^2$ is $\{0\}$ if $c=-1$ and $\{0,-1\}$ if $c=1$, so $S(X^4+cX^2+e)=e$ if $c=-1$ and $-1-e$ if $c=1$;
the $c=1$ case is consistent with Conjecture~\ref{conj3}, but the $c=-1$ case violates Conjecture~\ref{conj3}.
If $q=5$ then $S(X^4+cX^2+e) = 2+2c^2e$ if $c\ne 0$, all of which violate Conjecture~\ref{conj3}.

Magma verification:

K:=GF(3); {x^4-x^2:x in K} eq {0} and {x^4+x^2:x in K} eq {0,-1} and {&+[i:i in {x^4-x^2+e:x in K}] eq e:e in K} join {&+[i:i in {x^4+x^2+e:x in K}] eq -1-e:e in K} eq {true};
K:=GF(5); {&+[i:i in {x^4+c*x^2+e:x in K}] eq 2+2*c^2*e:c,e in K|c ne 0} eq {true};

\end{comment}
%

\begin{proof}[Proof of Corollary~\ref{conj2}]
Plainly if $a=0$ then the image of $\F_q$ under $X^4+aX^2$ consists of $0$ and the set $W$ of fourth powers in $\F_q^*$, so that $S(X^4+aX^2)=0$.  If $q\equiv 1 \pmod 8$ then by Theorem~\ref{conj3} we have
\[
T = \Bigl\{\frac{-9a^2}{64} : a\in\F_q,\, \chi(a)\in\{0,1\}\Bigr\} \cup \Bigl\{\frac{-a^2}{64} : a\in\F_q,\,\chi(a)=-1\Bigr\}.
\]
Since $-1$ and $64$ are fourth powers in $\F_q$, this equals 
\[
\{9x^4 : x \in \F_q\} \cup \{x : x \in \F_q\setminus W,\, \chi(x)=1\}.
\]
If $\charp(\F_q)=3$ then $T$ consists of the squares in $\F_q$ which are not in $W$.  By quadratic reciprocity, if $\charp(\F_q)\ne 3$ then $\chi(3)=1$ if and only if $q\equiv 1\pmod 3$, so if $q\equiv 1 \pmod 3$ then $T$ consists of the squares in $\F_q$, and if $q\equiv 2 \pmod 3$ then $T$ consists of the squares in $\F_q$ which are not in $W$.

If $q\equiv 5 \pmod 8$ then 
\[
T = \Bigl\{\frac{-a^2}{64} : \chi(a)\in\{0,1\}\Bigr\} \cup \Bigl\{\frac{-9a^2}{64} : \chi(a)=-1\Bigr\}.
\]
Here $-1$ and $64$ are squares which are not fourth powers, so 
\[
T = \{0\} \cup W \cup \{9x : x\in\F_q\setminus W,\, \chi(x)=1\}.
\]
If $\charp(\F_q)=3$ then $T=\{0\}\cup W$.  If $\charp(\F_q)\ne 3$ then $\chi(3)=1$ if and only if $q\equiv 1\pmod 3$, so if $q\equiv 1 \pmod 3$ then $T$ consists of the squares in $\F_q$, and if $q\equiv 2 \pmod 3$ then $T$ consists of the fourth powers in $\F_q$.

If $q\equiv 3 \pmod 4$ then $T = \{\epsilon x : \chi(x)\in\{0,1\},\, \epsilon \in \{1,-7\}\}$.  By quadratic reciprocity, $-7$ is a square in $\F_q$ if and only if $q$ is a square in $\F_7$.
Thus if $q\equiv 0,1,2,4 \pmod 7$ then $T$ consists of the squares in $\F_q$, and if $q\equiv 3,5,6 \pmod 7$ then $T=\F_q$.  Finally, if $a=8$ then $\chi(-2a)=-1$, so $S(X^4+8X^2)=8^2/64=1$.
\end{proof}

%#######################################################################
%#######################################################################

\section{Proofs of remaining results}\label{sec4}

In this section we prove Propositions~\ref{3}, \ref{power}, and \ref{odd}.  We prove these in reverse order.

\begin{proof}[Proof of Proposition~\ref{odd}]
Let $\zeta\in\F_q^*$ have order $s$.  Then $f(\zeta X)=\zeta^r f(X)$, so that $f(\F_q)$ is preserved by multiplication by $\zeta^r$.  Since $s\nmid r$, we have $\zeta^r\ne 1$ and thus $S(f)=0$.
\end{proof}

Note that Corollary~\ref{oddplus} follows immediately from Proposition~\ref{odd}.

\begin{proof}[Proof of Proposition~\ref{power}]
The set $Z$ of $n$-th powers in $\F_q$ equals the set of $m$-th powers, where $m:=\gcd(n,q-1)$.  Thus $Z$ consists of $0$ and the $(q-1)/m$-th roots of unity, and plainly $f(\F_q)=\{ax+b:x\in Z\}$.
If $m\ne q-1$ then $Z$ contains an element $\delta\notin\{0,1\}$, so since $Z$ is preserved by multiplication by $\delta$ we conclude that the sum of the elements of $Z$ equals $0$, whence $S(f)=b\cdot\abs{Z}=b\cdot (q+m-1)/m=b\cdot (m-1)/m$.  Finally, if $m=q-1$ then $Z=\{0,1\}$ so that $f(\F_q)=\{b,a+b\}$ and thus $S(f)=a+2b$.
\end{proof}

\begin{lemma}\label{deg3vals}
Let $q$ be a prime power, and let $f(X):=X^3+aX$ with $a\in\F_q^*$.  If $q>3$ then
$\abs{f(\F_q)}=\lfloor(2q+1)/3\rfloor$.
\end{lemma}

Lemma~\ref{deg3vals} was first proved in full generality in \cite[Prop.~4.6]{T-val}; see \cite[Rmk.~4.8(b)]{T-val} for a discussion of earlier literature.

\begin{proof}[Proof of Proposition~\ref{3}]
If $f(X)=aX+b$ then $f(\F_q)=\F_q$, so that $S(f)$ 
% equals $0$ or $1$ according as $q>2$ or $q=2$.
equals $0$.

If $q$ is odd and $f(X)=aX^2+bX+c$  then $f(X)=a(X+b/(2a))^2+e$ where $e:=c-b^2/(4a)$, so that $f(\F_q)=\{ax^2+e:x\in\F_q\}$.  Since there are $(q+1)/2$ squares in $\F_q$,
% and their sum is $0$ if $q>3$ and $1$ if $q=3$, it follows that $S(f)=e/2$ if $q>3$ and $a-e$ if $q=3$.
and their sum is $0$, it follows that $S(f)=e/2$.

If $q$ is even and $f(X)=aX^2+c$  then $f(\F_q)=\F_q$, so that $S(f)$ equals 
% $0$ if $q>2$ and $1$ if $q=2$.
$0$.

Now suppose $q$ is even and $f(X)=aX^2+bX+c$ with $b\ne 0$.  Then $(a/b^2)\cdot f(bX/a) = g(X)+e$ where $g(X):=X^2+X$ and $e:=ac/b^2$.  
Plainly $g(X)$ induces a homomorphism on the additive group of $\F_q$ with kernel $\F_2$, so that $\abs{g(\F_q)}=q/2$.  It is well-known that $g(\F_q)$ consists of the roots of the squarefree polynomial $T(X):=X^{q/2}+X^{q/4}+\dots+X$.
Thus $S(g)$ is the coefficient of the term of $T(X)$ of degree $(q-2)/2$,
% so that $S(g)=1$ if $q=4$ and $S(g)=0$ if otherwise.  Finally, $(a/b^2)\cdot S(f)=S(g)+e\cdot\abs{S(g)}=S(g)+e\cdot q/2$, so that $S(f)=0$ if $q>4$, $S(f)=b^2/a$ if $q=4$, and $S(f)=b^2e/a=c$ if $q=2$.
so that $S(g)=0$.
Finally, $(a/b^2)\cdot S(f)=S(g)+e\cdot\abs{S(g)}=e\cdot q/2=0$.

Next suppose $3\nmid q$ and $f(X)=aX^3+bX^2+cX+d$.  Then $f(X-b/(3a)) = g(X)+e$ where $g(X):=aX^3+\tilde{c}X$ with $\tilde{c}:=(3ac - b^2)/(3a)$ and $e:=(27a^2d - 9abc + 2b^3)/(27a^2)$.
We may assume $\tilde{c}\ne 0$, since otherwise the result follows from Proposition~\ref{power}.
Thus $\abs{g(\F_q)}=\lfloor(2q+1)/3\rfloor$ by Lemma~\ref{deg3vals}, 
% {\color{blue}(The above assertion on $\abs{g(\F_q)}$ in general is true only if $q>3$, although the remaining case $q=2$ is trivial.)}
so since $S(f)=S(g)+e\cdot\abs{g(\F_q)}$ it remains to determine $S(g)$.  If $q$ is odd then $S(g)=0$ by Proposition~\ref{odd}.
Now suppose $q$ is even.  Then $g'(X)=aX^2+\tilde{c}=g(X)/X$, so that no element of $\F_q^*$ has exactly two $g$-preimages in $\F_q$.  Hence every element of $f(\F_q)\setminus\{0\}$ has an odd number of $\F_q$-preimages, so that $S(g)=\sum_{x\in\F_q} g(x)=0$.
% which equals $0$ when $q>4$, equals $a$ when $q=4$, and equals $a+\tilde{c}$ when $q=2$.

Now suppose $3\mid q$ and $f(X)=aX^3+cX+d$, and write $g(X):=aX^3+cX$.  By Corollary~\ref{oddplus} we have $S(f)=d\cdot \abs{g(\F_q)}$.
Here $g(X)$ induces a homomorphism on the additive group of $\F_q$, so that $\abs{g(\F_q)}=q/m$ where $m$ is the number of roots of $g(X)$ in $\F_q$.  
Plainly $m=3$ if $\chi(-c/a)=1$ and $m=1$ otherwise.

Finally, suppose $3\mid q$ and $f(X)=aX^3+bX^2+cX+d$ with $b\ne 0$.  Then $f(X+c/b)=g(X)+e$ where $g(X):=aX^3+bX^2$ and $e:=(ac^3+b^3d-b^2c^2)/b^3$.  We have $S(f)=S(g)+e\cdot\abs{g(\F_q)}$.
Since $g((b/a)\cdot X)=(b^3/a^2)\cdot h(X)$ where $h(X):=X^3+X^2$, we have $S(g)=(b^3/a^2)\cdot S(h)$ and $\abs{g(\F_q)}=\abs{h(\F_q)}$.  To conclude the proof, it remains to determine $S(h)$ and $\abs{h(\F_q)}$, which is done in Lemma~\ref{char3} below.
\end{proof}

\begin{lemma}\label{char3}
Let $q=3^k$ and $h(X):= X^3+X^2$. Then an element $\gamma\in\F_q$ has more than one $h$-preimage in\/ $\F_q$ if and only if $\gamma=\delta^2$ for some $\delta\in\F_q$ with\/ $\Tr_{\F_q/\F_3}(\delta)=0$, and $\gamma$ has three such $h$-preimages if and only if $\delta\ne 0$.  It follows that $\abs{h(\F_q)}=2q/3$, and that $S(h)=0$ if $k>2$ and $S(h)=-1$ if $k\le 2$. 
\end{lemma}

\begin{proof}
Since the only root of $h'(X)=2X$ is $0$, the only $\gamma\in\F_q$ with fewer than three distinct $h$-preimages in $\mybar\F_q$ is $h(0)=0$, whose $h$-preimages are $0$ and $-1$.  Writing $Z_i$ for the set of elements in $\F_q$ having exactly $i$ distinct $h$-preimages in $\F_q$, it follows that $Z_2=\{0\}$.
We claim that $Z_3$ consists of the elements $\delta^2$ where $\delta\in\F_q^*$ satisfies $\Tr_{\F_q/\F_3}(\delta)=0$.  We now deduce the value of $\abs{h(\F_q)}$ from this claim.  The claim implies that $\abs{Z_3}=(q-3)/6$, and that $\abs{h^{-1}(Z_3)}=(q-3)/2$. Thus
\[
\abs{Z_1}=\abs{\F_q\cap h^{-1}(Z_1)}=\abs{\F_q\setminus h^{-1}(Z_2\cup Z_3)}=q-2-\frac{q-3}{2} = \frac{q-1}{2},
\]
so that
\[
\abs{h(\F_q)} = \abs{Z_1}+\abs{Z_2}+\abs{Z_3}=\frac{q-3}{2}+1+\frac{q-1}{6} = \frac{2q}{3}.
\]

Now we prove the claim. For $\gamma\in\F_q^*$, we have $\gamma\in Z_3$ if and only if $X^3+X^2-\gamma$ has three roots in $\F_q$, or equivalently $\gamma X^3-X-1$ has three roots in $\F_q$. Since $\gamma X^3-X$ induces a homomorphism on the additive group of $\F_q$, the polynomial $\gamma X^3-X-1$ has three roots in $\F_q$ if and only if it has at least one such root and $\gamma X^3-X$ has three roots in $\F_q$.  
The last condition holds if and only if $\gamma=\delta^2$ for some $\delta\in\F_q^*$.  Upon substituting $X/\delta$ for $X$, we see that $\delta^2 X^3-X-1$ has a root in $\F_q$ if and only if $X^3-X-\delta$ has such a root, or equivalently $\Tr_{\F_q/\F_3}(\delta)=0$.

It remains to determine $S(h)$.  This is easy when $k\le 2$, so we now assume $k>2$.  
Write $T(X):=X^{q/3}+X^{q/9}+\dots+X$, so that $T(X)$ is squarefree with roots being the elements $\delta\in\F_q$ for which $\Tr_{\F_q/\F_3}(\delta)=0$.
Then $T(X)=X R(X^2)$ where $R(X):=X^{(q/3-1)/2}+X^{(q/9-1)/2}+\dots+1$, so the sum of the squares of the roots of $T(X)$ equals the sum of the roots of $R(X)$, which is zero since $R(X)$ has no term of degree $\deg(R)-1$.  Hence the sum of the elements of $Z_3$ is zero, so since $Z_2=\{0\}$ and $\sum_{x\in\F_q} h(x)=0$ we conclude that $S(h)=0$.
\end{proof}

\begin{rmk}
The value of $\abs{h(\F_q)}$ in Lemma~\ref{char3} was first determined in \cite{CM} and \cite[Prop.~4.7]{T-val}.  In fact, \cite{CM} determines $\abs{g(\F_q)}$ whenever $g(X)=X^{Q-1}(X+1)$ with $q=Q^k$, and one can combine further results from \cite{CM} with arguments in the above proof to show that $S(g)=0$ when $k>2$ or $q=2$, and $S(g)=-1$ otherwise.
\end{rmk}

\end{document}